\documentclass [12pt]{amsart}

\usepackage{amssymb,amsxtra,amsfonts}
\usepackage{amsmath, accents} %
\usepackage{graphics}

\usepackage[colorlinks]{hyperref}

\usepackage{xspace} 

\long\def\pb #1*/{}

\def\reE@DeclareMathSymbol#1#2#3#4{%
    \let#1=\undefined
    \DeclareMathSymbol{#1}{#2}{#3}{#4}}
\DeclareSymbolFont{symbolsC}{U}{txsyc}{m}{n}
\SetSymbolFont{symbolsC}{bold}{U}{txsyc}{bx}{n}
\DeclareFontSubstitution{U}{txsyc}{m}{n}
\reE@DeclareMathSymbol{\strictiff}{\mathrel}{symbolsC}{76}

\newcommand\beq{\begin{equation}}
\newcommand\eeq{\end{equation}}
\newcommand\bal{\begin{align*}}
\newcommand\eal{\end{align*}}   
\newcommand\bmx{\left(\begin{matrix}}
\newcommand\emx{\end{matrix}\right)}
\newcommand\bsmx{\left(\begin{smallmatrix}}
\newcommand\esmx{\end{smallmatrix}\right)}
\newcommand\bmxnp{\begin{matrix}}
\newcommand\emxnp{\end{matrix}}
\newcommand\bsmxnp{\begin{smallmatrix}}
\newcommand\esmxnp{\end{smallmatrix}}

\newcommand{\wt}{\widetilde}

\DeclareMathSymbol{\widehatsym}{\mathord}{largesymbols}{"62}

\renewcommand{\d}{\partial}

\newcommand{\bSi}{{\bf \Si}}

\newcommand{\onto}{\twoheadrightarrow}

\newcommand{\into}{\hookrightarrow}
\newcommand{\spq}{/\!\!/}

\newcommand{\st}{\ \bigl\vert\ }

\newcommand{\range}{\operatorname{Im}}

\providecommand{\<}{\langle}
\renewcommand{\>}{\rangle}

\def\part#1{\frac{\partial\phantom{q}}{\partial#1}}

\newcommand{\fus}{\circledast}
\newcommand{\fusion}[1]{\underset{#1}{\circledast}}

\DeclareFontFamily{U}{wncy}{}
\DeclareFontShape{U}{wncy}{m}{n}{<->wncyr10}{}
\DeclareSymbolFont{mcy}{U}{wncy}{m}{n}
\DeclareMathSymbol{\Sh}{\mathord}{mcy}{"58} 

\usepackage{xparse}
\NewDocumentCommand{\dn}{e{_^}}{
  _{\IfValueT{#1}{#1}\vphantom{\smash[b]{|}}}
  ^{\IfValueT{#2}{#2}\vphantom{\smash[t]{|}}} }

\newcommand{\MB}{\mathcal{M}_{\text{\rm B}}}

\newcommand{\Id}{\text{\rm Id}}

\newcommand{\Gal}{{\rm Gal}} 

\newcommand{\Ad}{{\mathop{\rm Ad}}}

\newcommand{\Ram}{\mathop{\rm Ram}} 

\newcommand{\slope}{\mathop{\rm slope}}

\newcommand{\Prod}{\prod}

\DeclareMathOperator{\Perm}{Perm}

\DeclareMathOperator{\Hom}{Hom}         
\DeclareMathOperator{\Aut}{\mathop{\rm Aut}}

\DeclareMathOperator{\GrAut}{GrAut}

\DeclareMathOperator{\THom}{THom}

\newcommand{\SL}{{\mathop{\rm SL}}}

\newcommand{\GL}{{\mathop{\rm GL}}}

\newcommand{\Out}{\mathop{\rm Out}}
\DeclareMathOperator{\Inn}{Inn} 

\newcommand{\Ker}{\mathop{\rm Ker}}

\newcommand{\ba}{{\bf a}}

\newcommand{\bC}{{\bf C}}

\newcommand{\bD}{{\bf D}}

\newcommand{\bk}{{\bf k}}

\newcommand{\bH}{{\bf H}}\newcommand{\bh}{{\bf h}}

\newcommand{\bK}{{\text{\bf K}}}

\newcommand{\bM}{{\bf M}}

\newcommand{\bt}{{\bf t}}
\newcommand{\bT}{{\bf T}}

\newcommand{\bx}{{\bf x}}
\newcommand{\bfx}{{\bf x}}
\newcommand{\bfy}{{\bf y}}

\DeclareSymbolFont{bbold}{U}{bbold}{m}{n}
\DeclareSymbolFontAlphabet{\mathbbold}{bbold}

\newcommand{\IA}{\mathbb{A}}

\newcommand{\IC}{\mathbb{C}}
\newcommand{\ID}{\mathbb{D}}

\newcommand{\IH}{\mathbb{H}}

\newcommand{\IL}{\mathbb{L}}
\newcommand{\IM}{\mathbb{M}}

\newcommand{\IP}{\mathbb{P}}

\newcommand{\IS}{\mathbb{S}}

\newcommand{\IT}{\mathbb{T}}

\newcommand{\IV}{\mathbb{V}}

\newcommand{\IZ}{\mathbb{Z}}

\newcommand{\cA}{\mathcal{A}}

\newcommand{\cC}{\mathcal{C}}

\newcommand{\cE}{\mathcal{E}}

\newcommand{\cF}{\mathcal{F}}
\newcommand{\cG}{\mathcal{G}}

\newcommand{\cI}{\mathcal{I}}

\providecommand{\cL}{}
\renewcommand{\cL}{\mathcal{L}}

\newcommand{\cP}{\mathcal{P}}

\providecommand{\cR}{}
\renewcommand{\cR}{\mathcal{R}}

\newcommand{\cT}{\mathcal{T}}

\usepackage[mathscr]{euscript}

\newcommand{\g}{       \mathfrak{g}     }

\newcommand{\al}{\alpha}

\newcommand{\be}{\beta}
\newcommand{\ga}{\gamma}

\newcommand{\Ga}{\Gamma}

\newcommand{\La}{\Lambda}

\newcommand{\Si}{\Sigma}
\newcommand{\Th}{\Theta}

\newcommand{\ov}{\overline}

\renewcommand{\bar}{\overline}

\makeatletter
 \newlength{\typesize}
\makeatother

\newlength{\vvoff}
\newlength{\hhoff}

\def\underset#1#2{\ \smash{\mathop{ #2 }\limits_{#1}}\ }

\newcommand{\pf}{\begin{bpf}}

\newcommand{\pfms}{\begin{bpfms}}
\newcommand{\epf}{\end{bpf}\hfill$\square$\\}           
\newcommand{\epfms}{\end{bpfms}\hfill$\square$\\}       

\newcommand{\idea}{\begin{bidea}}

\newcommand{\eidea}{\end{bidea}\hfill$\square$\\}           

\newcommand{\sk}{\begin{bsk}}    

\newcommand{\esk}{\end{bsk}\hfill$\square$\\}           
\newcommand{\sketch}{\begin{bsketch}}

\newcommand{\esketch}{\end{bsketch}\hfill$\square$\\}

\newtheorem {hypo}{\bf\hspace{-\parindent}Hypothesis}
\newtheorem {thm}[hypo]{Theorem}   
\newtheorem {prop}[hypo]{Proposition}

\newtheorem {cor}[hypo]{Corollary}
\newtheorem {lem}[hypo]{Lemma}

\theoremstyle{definition}\newtheorem {defn}[hypo]{Definition}
\theoremstyle{definition}
\theoremstyle{definition}\newtheorem{eg}[hypo]{Example} 
\theoremstyle{remark}\newtheorem{rmk}[hypo]{Remark}
\theoremstyle{remark}

\renewcommand{\ba}{\al}  %

\begin{document}

\date{December 3, 2025}

 \title[Polystability of Stokes representations and Galois groups]{Polystability of Stokes representations and differential Galois groups}
 \author{Philip Boalch and Daisuke Yamakawa}%

\begin{abstract}
Polystability of (twisted) Stokes representations (i.e. wild monodromy representations) will be characterised, 
in terms of the corresponding differential Galois group 
(generalising the Zariski closure of the monodromy group
in the tame case).
This  extends some results of Richardson. 
Further, the intrinsic approach to such results will be established, 
in terms of reductions of Stokes local systems.
\end{abstract}

 \maketitle

\setcounter{tocdepth}{1}
\tableofcontents

\section{Introduction}

We continue our investigations of the nonabelian moduli spaces 
in 2d gauge theory (the theory of connections on curves). 
This article is a continuation of the paper \cite{twcv} 
that completed the construction of the wild character varieties 
(irregular Betti spaces) of smooth curves 
as affine algebraic Poisson varieties 
(using the ``quasi-Hamiltonian'' TQFT approach to meromorphic connections, completing the sequence \cite{saqh, fission, gbs} that started in 2002). 
The scope of \cite{twcv} encompasses the monodromy/Stokes data 
of meromorphic connections with arbitrary formal structure at the poles 
(and arbitrary, possibly twisted, complex reductive structure group). 
This is a purely algebro-geometric approach, 
complementary to the earlier analytic approaches 
\cite{wnabh, smid} 
that constructed wild nonabelian Hodge moduli spaces as 
(complete) hyperk\"ahler manifolds, 
and identified these moduli spaces of wild harmonic bundles 
as both spaces of meromorphic connections 
and meromorphic Higgs bundles on parabolic vector 
bundles %
(encompassing  
the wild topological symplectic structures on the De\! Rham/Betti side  \cite{thesis, smid, saqh},  
the [good symplectic leaves of the] 
Bottacin--Markman algebraic Poisson structures on the 
meromorphic Higgs side \cite{Bot, Mar}, 
and the enriched Riemann--Hilbert--Birkhoff correspondence
of \cite{smid} Corollary 4.9.)
Here we will give an intrinsic characterisation of 
the points of the wild character varieties 
(i.e.\ the polystable Stokes representations), 
generalising results previously only available in the tame case, 
and we will also characterise the stable points 
(generalising a result of \cite{gbs} in the untwisted wild case). 
As in Richardson's paper \cite{Richardson} the whole paper is, 
to some extent, an exercise in the Hilbert--Mumford theorem 
(although, somehow, it took us an incredibly long time to 
correctly formulate the main intrinsic statements below).  
For more background and applications see the reviews in 
\cite{hdr,ihptalk, p12,hit70, tops} 
(the first large class of examples of wild character varieties 
is due to Birkhoff \cite{birkhoff-1913}, clarified in \cite{BJL79,JMU81,smid}, 
and the simplest case underlies the Drinfel'\!d--Jimbo quantum group $U_q (\g)$, 
as reviewed in \cite{hdr} \S4).

In particular the reader should bear in mind the important results 
of Deligne \cite{Del70} telling us, in particular, 
that the complex representations of the topological fundamental group 
of a smooth complex algebraic curve 
parameterise purely algebraic objects known as 
``algebraic connections with regular singularities at the punctures'' 
(see also the reviews \cite{BoreletalDmods} Ch.~4 (Thm.~7.2.1), 
\cite{Katz76} by Malgrange and Katz). 
The theory of wild character varieties stems from the deep work (1857-1991), 
summarised in \cite{tops}, that gives the {\em generalisation} 
of the  topological fundamental group, 
whose representations parameterise all the other algebraic connections, 
i.e.\ those with {\em irregular} singularities, 
having a plethora of applications. 
The simplest point of view seems to be the wild surface groupoids \cite{gbs, twcv}  
(also reviewed in \cite{tops,c236}) which generalise 
the Poincar\'e fundamental groupoids (with a finite number of basepoints) 
and easily yield the wild character varieties 
(i.e.\ the irregular Betti moduli spaces) 
and the wild character stacks (\cite{tops} \S13.2, i.e.\ the irregular Betti stacks).
Note also that it is the consideration of (poly)stability conditions (not the stack) 
that leads to the algebraic parametrisation of wild harmonic bundles, 
following in line with a long sequence of generalisations starting with the  
Narasimhan--Seshadri and Mehta--Seshadri theorems 
(giving purely algebraic approaches to unitary flat $C^\infty$ connections, 
i.e.\ unitary representations of the fundamental group) 
\cite{NarSes, MehSes}.

First we will recall the basic statements in the tame case.
Let $G$ be a connected complex reductive group, such as $\GL_n(\IC)$, and let $\Si^\circ$ be a smooth complex algebraic curve.
Thus $\Si^\circ=\Si\setminus \ba$ 
for some smooth compact complex algebraic curve $\Si$ (i.e. a compact Riemann surface), and a finite subset $\ba\subset \Si$. 
Given a basepoint $b\in \Si^\circ$ one can consider the representation variety
$$\cR=\Hom(\pi_1(\Si^\circ,b),G)$$
which is a complex affine variety equipped with an action of $G$, conjugating representations.
A point $\rho\in \cR$ is a $G$-representation, i.e. a group homomorphism 
$\rho:\pi_1(\Si^\circ,b)\to G$. 
In turn the $G$-character variety, or Betti moduli space
$$\MB=\MB(\Si^\circ,G)=\cR/G=\Hom(\pi_1(\Si^\circ,b),G)/G$$
is the (affine) geometric invariant theory quotient of $\cR$ by $G$.
By definition this means that the points of the variety 
$\MB$ are the {\em closed} $G$ orbits in $\cR$.
The representations $\rho$ whose  $G$ orbits are closed are called the 
{\em polystable} representations.
A basic question is thus to characterise the polystable representations intrinsically.
The answer (due to Richardson, building on earlier work in the general linear case), is as follows.

Let $\Gal(\rho)=A(\rho)\subset G$ be the 
Zariski closure ({\em adh\'erence}) of the image of $\rho$.
A theorem of Schlesinger \cite{schleshbk21} p.101 
implies that if $\rho$ is the monodromy representation of an algebraic connection $(\nabla,E)$ on an algebraic principal $G$-bundle $E\to \Si^\circ$, with regular singularities at $\ba$, then 
$\Gal(\rho)$ is the differential Galois group of $(\nabla,E)$.
Recall %
that a complex affine algebraic group is a {\em linearly reductive group} if its identity component is reductive (i.e. has trivial unipotent radical).
The basic characterisation is then:

\begin{thm}[\cite{Richardson}]
A point $\rho\in \cR$ is polystable if and only if $\Gal(\rho)$ is a linearly reductive group.
Further, $\rho$ is stable if and only if $\Gal(\rho)$ is not contained in a proper parabolic subgroup of $G$. 
\end{thm}

A representation $\rho\in \cR$ whose Galois group  $\Gal(\rho)$ is linearly reductive is often called a {\em semisimple} 
representation\footnote{Indeed if $V\cong \IC^n$ and $G=\GL(V)\cong\GL_n(\IC)$ then 
$\rho:\pi_1(\Si^\circ,b)\to G$ is a {semisimple} representation
if and only if $V$ is a semisimple 
$\pi_1(\Si^\circ,b)$-module.}.
Thus the theorem says that $\rho\in \cR$ is polystable if and only if 
$\rho$ is semisimple, and so the points of  $\MB$ are the isomorphism classes of semisimple representations.
Note that for $G=\GL_n(\IC)$ this statement is older (Artin/Procesi) and a full account appears in the book of Lubotzky--Magid \cite{lubmagid}.

This paper is concerned with the extension of this result to the case of Stokes representations, generalising the fundamental group representations, and the 
characterisation of the 
points of the wild character varieties, 
generalising the (tame) character varieties appearing above.

In brief, any algebraic connection $(\nabla,E)$ on an algebraic  principal $G$-bundle $E\to \Si^\circ$, 
has an invariant, its irregular class, at each point 
$a\in \ba$. A connection is regular singular if and only if its irregular class is trivial.
Our aim is to give the generalisation of Richardson's results when the irregular classes are arbitrary.
Due to work of many people 
it is known how to describe algebraic 
connections completely 
topologically, so we can work algebraically on the Betti/Stokes side
(see the review of this story 
in \cite{tops}).
The basic notions from the tame case are generalised as follows:
\begin{align*}
\text{local system} \ &\rightsquigarrow \ 
\text{Stokes local system}\\
\text{ fundamental group (or groupoid)} 
 \ &\rightsquigarrow \  
 \text{wild surface group (or groupoid)}\\
\text{fundamental group representation } 
\ &\rightsquigarrow \ 
 \text{Stokes representation.}
 \end{align*}

Once these generalisations are understood then the story proceeds  similarly to the tame case
(defining a wild representation variety $\cR$ parameterising 
framed Stokes local systems, 
and then acting by a reductive group 
to forget the framings).
A key novelty in the wild setting is that there is a breaking of structure group (``fission'') near the 
marked points so the group acting involves a  (reductive) subgroup of $G$. 
This is intimately related to extra generators in Ramis' description of the differential Galois group
\cite{ramis-pv}, 
generalising Schlesinger's density theorem. 
However, as we will recall, 
the basic feature of the tame case remains, that, 
by choosing a presentation of the wild surface groupoid, 
the wild representation variety 
$\cR$ can be written explicitly in terms of a product 
of simple explicit  pieces (internally fused doubles $\ID$ 
or fission spaces $\cA$), 
one for each marked point or handle on $\Si$:
\beq
\cR \ \cong\  (\ID^{\fus{}g}\fus{}\cA_1\fus\cdots\fus\cA_m )\spq G\eeq
if $\Si$ has genus $g$ and $m$ marked points, where $\fus{}$
is the quasi-Hamiltonian fusion product.

\subsection{Summary of main results}

In this section we will summarise the main results in a short uncluttered form, with links to the references to the full definitions (mainly in \cite{twcv}). 
More details are in the body of the article, and 
an overview of the simpler $\GL_n(\IC)$ set-up is in \cite{tops} \S13.

Let $\Si$ be a smooth compact complex algebraic curve (i.e. a compact Riemann surface) and 
let $\ba\subset \Si$ be a finite non-empty subset. 
Let $G$ be a connected complex reductive group.

Choose a $G$-irregular class $\Th_a$ at 
each point $a\in \ba$ 
(possibly twisted) as in \cite{twcv} \S3.5.  
Let $\bSi=(\Si,\ba,\Th)$ be the 
resulting wild Riemann surface 
with structure group $G$ as in \cite{twcv} \S4 (and \cite{gbs} \S8.1, Rmk 10.6).   
The wild Riemann surface $\bSi$ 
has a fundamental groupoid $\Pi$ and a 
 character variety $\MB(\bSi)$, as follows.
 
First the topological 
notion of Stokes local 
system $\IL$ on $\bSi$ is defined as in 
\cite{twcv} Defn. 13 (and \cite{gbs} Rmk. 8.4, \cite{p12} \S4.3, \cite{tops} \S8).

As in \cite{twcv} \S5, \cite{gbs} \S8, it follows that 
if we choose suitable basepoints $\be$  then the
wild surface groupoid $\Pi=\Pi_1(\bSi,\be):=\Pi_1(\wt\Si,\be)$ 
is well defined,   as is the
wild representation variety (the space of Stokes representations of $\Pi$):
$$\cR(\bSi,\be) = \Hom_\IS(\Pi,G)\subset \Hom(\Pi,G).$$ 
It is an affine variety equipped with an action of a complex reductive group $\bH$. 

Any framed Stokes local system $\IL$ determines a Stokes representation  $\rho=\rho_\IL\in\cR(\bSi,\be)$. 
Two Stokes local systems are isomorphic if and only if their Stokes representations are in the same 
$\bH$-orbit in $\cR(\bSi,\be)$.

The notion of $\rho$ being polystable or stable for the action of $\bH$ on $\cR$ is well-defined, as for any action of a reductive group on an affine variety (as in \cite{Richardson}).
The points of the Poisson wild character variety $\MB(\bSi)$ are the polystable (i.e. closed) $\bH$-orbits.

On the other hand $\rho$ determines the Galois group
$\Gal(\rho)\subset G$, as in Ramis' density theorem 
\cite{ramis-pv}
(see \S\S\ref{ssn: gg from sreps}, \ref{ss: galois group} below).
It involves not just the usual monodromy, 
but also the formal monodromy, Stokes automorphisms and the Ramis  tori.

Finally we can define $\IL$ to be 
irreducible if it has 
no proper parabolic reductions, and to be reductive (or ``semisimple'') if it has an irreducible Levi reduction.

\begin{thm}\label{thm: 1to3}
Let $\IL$ be a Stokes local system on $\bSi$, and let 
$\rho$ be its Stokes representation.
The following are equivalent:

$\bullet 1)$\ $\rho$ is polystable,

$\bullet 2)$\ $\Gal(\rho)$ is linearly reductive,

$\bullet 3)$\ $\IL$ is a semisimple Stokes local system.
\end{thm}

Thus the points of the wild character variety  $\MB(\bSi)$
are the isomorphism classes of semisimple Stokes local systems on $\bSi$.
In the tame case (with each irregular class trivial) 
the groupoid $\Pi$ becomes Poincar\'e's fundamental groupoid (with a finite number of basepoints), and the theorem is already known \cite{Richardson}.

Further we will consider stability (not just polystability).
This requires possibly adding  a few extra punctures to 
control the kernel of the action, but no generality is lost 
(see the discussion after \eqref{eq: stokes gluing}).

\begin{thm} \label{thm: 4to6}
Let $\IL$ be a Stokes local system on $\bSi$, and let 
$\rho$ be its Stokes representation.
The following are equivalent:

$\bullet 4)$\ $\rho$ is a stable point of $\cR(\bSi,\be)$,

$\bullet 5)$\ $\Gal(\rho)$ is not contained in a proper parabolic subgroup of $G$,

$\bullet 6)$\ $\IL$ is irreducible.
\end{thm}

This result was already 
established in \cite{gbs} in the case
where each irregular class was not twisted.

Recall from \cite{twcv} that two types of twist are possible: the formal twists (twisted irregular classes, as above), and also interior twists, over the interior of the curve, 
where we start with a local system of groups 
$\cG\to \Si^\circ$, with each fibre isomorphic to $G$.
Similarly we will establish the analogous results 
in this fully twisted setting.
We will suppose that $\cG$ is ``out-finite'' 
in the sense of \S\ref{ss: grlocsys} below. 

The  first main difference 
(in the presence of interior twists) 
is that the wild representation variety 
is replaced by a space of twisted Stokes 
representations:
$$\cR=
\THom_\IS(\Pi,G)\subset 
\Hom(\Pi,G\ltimes \Aut(G))$$
which is an affine variety equipped with an action of a complex reductive group $\bH$.  
Secondly $\Gal(\rho)$ is not a subgroup of $G$, but rather it comes with a homomorphism 
$\Gal(\rho)\to \Aut(G)$, so naturally acts on $G$ by automorphisms. The image of $\Gal(\rho)$ in  $\Aut(G)$ 
will be denoted $\bar\Gal(\rho)$.
Then we can define 
irreducibility (\S\ref{ss irred gr loc sys}) and semisimplicity (\S\ref{ss: red-ss loc sys}) for 
Stokes $\cG$-local systems, and will 
prove analogues of the above results:

\begin{thm}\label{thm: primedversions}
Let $\bSi=(\Si,\ba,\Th)$ 
be a wild Riemann surface with 
group $\cG\to \Si^\circ$.
Let $\IL$ be a Stokes $\cG$-local system on $\bSi$, 
and let 
$\rho\in \THom_\IS(\Pi,G)$ be its twisted Stokes representation.
The following three conditions are equivalent:

$\bullet 1')$\ $\rho$ is a polystable point of 
$\cR(\bSi,\be)=\THom_\IS(\Pi,G)$,

$\bullet 2')$\ $\bar\Gal(\rho)$ is linearly reductive,

$\bullet 3')$\ $\IL$ is a semisimple Stokes $\cG$-local system.

\noindent Moreover the following  three conditions 
are also equivalent:

$\bullet 4')$\ $\rho$ is a stable point of 
$\cR(\bSi,\be)=\THom_\IS(\Pi,G)$,

$\bullet 5')$\ $\bar\Gal(\rho)\subset \Aut(G)$ does not preserve a proper parabolic subgroup of $G$,

$\bullet 6')$\ $\IL$ is an irreducible Stokes $\cG$-local system.
\end{thm}

In the set-up of Theorems \ref{thm: 1to3},\ref{thm: 4to6} with $\cG$ constant it is true that 
$a')\Leftrightarrow a)$ for all $a=1,2,\ldots 6$. 
Thus  Thm. \ref{thm: primedversions} 
implies both Theorems \ref{thm: 1to3},\ref{thm: 4to6}.

\subsection{Layout of the article}\

Sections \ref{sn: trich}, \ref{sn: moregen} generalise 
some of Richardson's results, in two steps. 
These are extrinsic results, and 
may well have further applications, beyond the wild character varieties of curves.

Section \ref{sn: apwcv} 
then applies these results to the 
spaces of Stokes representations leading to the 
wild character varieties $\MB(\bSi)=\cR/\bH$.
The main results are the equivalences 
$1')\Leftrightarrow 2')$ in Cor. \ref{cor: main-red-annonce}, and $4')\Leftrightarrow 5')$ in Cor. 
\ref{cor: irredGalrho}.

Section \ref{sn: SLS}
then discusses the intrinsic objects, 
Stokes local systems, and how stability/polystability 
can be read off in terms of (twisted) 
reductions of structure group.
The main results are the equivalences 
$1')\Leftrightarrow 3')$ and 
$4')\Leftrightarrow 6')$ in parts 2) and 
1) of Thm. \ref{thm: red-irred-SLS} respectively.

\subsection{Some other directions}
Note that for constant general linear groups the irreducible Stokes local systems are equivalent to the input data in the construction of wild harmonic metrics in \cite{sab99} (via \cite{Sim-hboncc}), 
so determine points of the wild nonabelian Hodge hyperk\"ahler manifolds \cite{wnabh}
(the full setting of \cite{wnabh} also allows non-zero Betti weights, as emphasised in \cite{wnabh} Rmk. 8.2).  
For other groups $G$ 
the untwisted irreducible Stokes local systems 
(i.e. as in \cite{gbs})
fit into the  ``Betti weight zero'' case of the 
$G$ wild nonabelian Hodge correspondence studied by Huang--Sun \cite{huangsun23}, which looks to be in line with our general conjecture in \cite{hit70} Rmk. 6, that the ``good'' meromorphic connections/Higgs fields on parahoric torsors are the right objects to look at.
This conjecture encompasses all the irreducible Stokes local systems 
considered here, and says further that 
the methods of \cite{wnabh} should yield  
complete hyperk\"ahler orbifolds as moduli spaces of such 
wild harmonic bundles, 
leading to  the 
full structure of ``nonabelian Hodge space'' as in \cite{wnabh} Defn. 7
(combining the key insights of both Hitchin and Simpson). 
For nontrivial Betti weights, one can also give an algebraic
construction of wild Betti spaces 
as in \cite{gbs} Rmk A.5, \cite{tops} \S 1.5 (1), \cite{yamakawa-mpa}. 

In another direction one of the  
key motivations for this work is the fact that
any admissible deformation 
of a wild Riemann surface leads 
to a local system of wild character varieties, 
and its monodromy
generalises the 
usual mapping class group actions on character varieties in the tame case.
A recent review of this story appears in \cite{twmcg}. 
As explained in \cite{smid, smapg}
the motivating examples for this whole 
line of thought 
were the Dubrovin--Ugaglia Poisson varieties 
whose 
braid group actions come from the braiding of 
counts of BPS states (\cite{Dub95long} Rmk 3.10, related to earlier work of Cecotti--Vafa); 
these are examples of {\em twisted} wild character varieties, involving the non-trivial outer automorphism of $\GL_n(\IC)$ (so we now have an
intrinsic general 
framework encompassing such examples).
Indeed it was by forgetting this twist that the 
Poisson variety $G^*$ underlying $U_q(\g)$ was 
recognised and identified 
as a framed wild character variety \cite{smapg}
(i.e. as a wild representation variety). 
See \cite{DRT22,DR23, DRY25} 
for some other recent developments 
concerning the generalised braid groups that act 
on wild character varieties, from admissible deformations of more general wild Riemann surfaces.

{\tiny 
{\bf Acknowledgements.} \parskip=0pt 
These results were completed in 2019 before the first named author moved departments and then 
learnt of the thesis work leading to 
\cite{cheng23}, that has some 
overlap with this paper in the tame setting with interior twists, 
although expressed in a slightly different language.  
In the intervening years we have not yet managed to 
incorporate possible simplifications suggested by \cite{cheng23} but thought it reasonable to release our original approach anyway since the scope is larger.

The second named author was supported by JSPS KAKENHI Grant Number 18K03256.

}

\section{Twisted version of Richardson's results}\label{sn: trich}

As in \cite{Richardson} we use the terminology that an 
affine algebraic group over $\IC$
is {\em reductive} if it is connected and has trivial unipotent radical. 
It is {\em linearly reductive} 
if its identity component is reductive.
We will also use the notation of \cite{Richardson} \S12
regarding semidirect products, such as $G\ltimes \Ga$ for 
$\Ga\subset\Aut(G)$.

Let $G$ be a linearly reductive group over $\IC$.
Recall that a point $x$ of an affine $G$-variety is said to be 
\emph{polystable} if the orbit $G \cdot x$ is closed. 
It is said to be \emph{stable} if it is polystable and  
the kernel of the action has finite index in its stabiliser $G_x$. 

Let $n$ be a positive integer.
For $\bfx =(x_i)_{i=1}^n \in G^n$, 
let $A(\bfx) \subset G$ be the Zariski closure of the subgroup 
generated by $x_1, x_2, \dots ,x_n$. 
In \cite{Richardson}, Richardson examined the simultaneous conjugation action of $G$ 
on the product $G^n$ and obtained the following results:

\begin{thm}[{\cite[Thm.~3.6]{Richardson}}]\label{thm:R1}
If $G$ is linearly reductive, a point $\bfx \in G^n$ is polystable   
if and only if $A(\bfx)$ is linearly reductive.
\end{thm}

\begin{thm}[{\cite[Thm.~4.1]{Richardson}}]\label{thm:R2}
If $G$ is reductive, a point $\bfx \in G^n$ is stable 
if and only if $A(\bfx)$ is not contained in any proper parabolic subgroup of $G$. 
\end{thm}

In this section we will establish twisted versions of these results, in Thm. \ref{thm:R1-tw} and Thm. \ref{thm:R2-tw} (2) respectively.
Two other characterisations of polystability will also be established, in Thm. \ref{thm:R2-tw} (1) and 
Cor. \ref{cor: completered}.
The subsequent section (\S\ref{sn: moregen}) will give a further generalisation.

Assume that $G$ is reductive and 
choose $\phi_1,\phi_2, \dots ,\phi_n \in \Aut(G)$. 
Let $\Gamma$ be the subgroup of $\Aut(G)$ 
generated by 
the inner automorphism group $\Inn(G)$ and 
$\phi_1,\phi_2, \dots ,\phi_n$. 
We assume that the quotient $\ov{\Gamma} := \Gamma/\Inn(G)$ is finite 
and regard $\Gamma$ as an algebraic group with 
identity component $\Inn(G) \cong G/Z(G)$.
For any $\phi\in \Ga$ let $G(\phi)$ denote the $G$-bitorsor/``twisted group'' 
$G\times\{\phi\}\subset G\ltimes \Ga$, as in 
\cite{twcv} \S2 (as explained there, the monodromy of a $\cG$-local system lies in such a space, and $G(\phi)$ embeds in the group of set-theoretic automorphisms of a fibre). 
Put 
\[
X = \prod_{i=1}^n G(\phi_i) \subset (G \ltimes \Gamma)^n,
\]
on which 
$G=G(\Id)$ 
acts by the simultaneous conjugation.
For $\bfx =(x_i)_{i=1}^n \in X$, 
let $A(\bfx)$ be the Zariski closure 
(in the algebraic group $G \ltimes \Gamma$) 
of the subgroup 
generated by $x_1, x_2, \dots ,x_n$. 
Let
$\ov{A}(\bfx) \subset \Gamma$ be the image of $A(\bfx)$ 
under the homomorphism
$$ \label{eq: Ad}
\Ad \colon G \ltimes \Gamma \to \Gamma, \quad x=(g,\phi) \mapsto 
\Ad(x) = \Ad(g) \circ \phi,
$$ 
where $\Ad(g) = g (\,\cdot\, )g^{-1}$. 
Note that $\ov{A}(\bfx)$ 
is contained in $\Aut (G)$ and 
hence naturally acts on $G$.
Thus in turn, via $\Ad$, the group 
$A(\bfx)$ naturally acts on $G$.

\begin{thm}\label{thm:R1-tw}
A point $\bfx \in X$ is polystable if and only if $\ov{A}(\bfx)$ is linearly reductive.
\end{thm}

\begin{rmk}
The notation $A(\bx)$ for the Zariski closure presumably stems from 
\cite{Borel56} \S3, where $A$ stands for {\em adh\'erence}.  
This may help to avoid possible confusion (since $\ov{A}(\bfx)$
denotes the image in $\Aut(G)$ here).
Note that \cite{Borel56} \S3.3 shows that if $H\subset G$ is any subgroup then 
the Zariski closure of the set $H$ is 
a Zariski closed subgroup of $G$. 
\end{rmk}

To see that Thm.~\ref{thm:R1-tw} generalises Thm.~\ref{thm:R1}, first note that:

\begin{prop}\label{prop: finiAAbar2}
Let $\wt G$ be a linearly reductive group with 
identity component $G$, and let $A$ be a closed subgroup of $\wt G$. 
Let $\ov A$ be the image of $A$ under the
map $\Ad\bigl\vert_G : \wt G \to\Aut(G)$, 
and regard it as a quotient algebraic group of $A$. Then $A$ is linearly reductive if and only if $\ov A$ is linearly reductive.
\end{prop}

\pf
Note that we may view $\Ad\bigl\vert_G$ as a map of algebraic groups: 
$\Aut(G)\cong \Inn(G)\ltimes \Out(G)$ may have an infinite number of components, but $\Ad\bigl\vert_G(\wt G)$ will only encounter a finite number of them.
Let $K$ be the kernel of the homomorphism 
$A \onto \ov{A}$. 
Note $K \cap G \subset Z(G)$,
and the identity component 
$K^0$ of $K$ is contained in $(\wt G )^0 = G$. 
Thus $K^0 \subset Z(G)$, 
which implies $K^0$ is a torus and so $K$ is linearly reductive. 
Now \cite{Richardson} 1.2.2 implies $A$ is linearly reductive if and only if $\ov A$ is linearly reductive.
\epf

\begin{lem}\label{lem: action cmpts}
Let $\wt{G}$ be a linearly reductive group with identity component $G$ 
and let $\bfx$ be a point in an affine $\wt{G}$-variety $Y$. 
Then the following hold.
\begin{enumerate}
\item $\bfx$ is polystable for the $\wt{G}$-action if and only if $\bfx$ is polystable 
for the $G$-action.
\item $\bfx$ is stable for the $\wt{G}$-action if and only if $\bfx$ is stable for the $G$-action.
\end{enumerate}
\end{lem}

\pf 
The $\wt{G}$-orbit $\wt{G} \cdot \bfx$ is a disjoint union
of a finite number of $G$-orbits, each of which is a $\wt{G}$-translate of $G \cdot \bfx$. 
Hence $\wt{G} \cdot \bfx$ is closed if and only if $G \cdot \bfx$ is closed.
The second assertion follows from 
the equality $G_\bfx = \wt{G}_\bfx \cap G$ for the stabilisers 
and a similar one for the kernels.
\epf

Thus in particular it is now clear that Thm.~ \ref{thm:R1-tw}  generalises Thm.~\ref{thm:R1}.

\pfms (of Theorem~\ref{thm:R1-tw}).
Since $\ov{\Gamma}$ is finite, we can find a finite subgroup 
$\Gamma'$ of $\Gamma$ such that $\Gamma = \Inn(G)\cdot \Gamma'$ 
thanks to a result of Borel--Serre and Brion 
(see \cite{borelserre64, brion15}).
For each $i=1,2, \dots ,n$ take $g_i \in G$ so that  
$\phi'_i := \Inn(g_i) \circ \phi_i \in \Gamma'$. 
Then we have isomorphisms of bitorsors
\[
f_i \colon G(\phi'_i) \to G(\phi_i), \quad (g,\phi'_i) \mapsto (g g_i, \phi_i)
\quad (i=1,2, \dots ,n),
\]
which induce a $G$-equivariant isomorphism
\[
f \colon X':= \prod_{i=1}^n G(\phi'_i) \to X, 
\quad (x'_i)_{i=1}^n \mapsto (f_i(x'_i))_{i=1}^n.
\]
Since $f$ is equivariant, a point $\bfx \in X$ is stable (resp.\ polystable) 
if and only if $f^{-1}(\bfx) \in X'$ is stable (resp.\ polystable).
Also, observe that $\Ad \circ f_i = \Ad$ ($i=1,2, \dots ,n$) and hence 
$\ov{A}(\bfx) = \ov{A}(f^{-1}(\bfx))$ for all $\bfx \in X$.
Therefore without loss of generality  
we may assume that $\phi_i \in \Gamma'$ for all $i=1,2, \dots ,n$.

Put $\wt{G} = G \ltimes \Gamma'$. It is a linearly reductive group 
with identity component $G$ and $X$ is a closed subvariety of $\wt{G}^n$. 
By Theorem~\ref{thm:R1}, a point $\bfx \in \wt{G}^n$ is polystable 
with respect to the 
simultaneous $\wt{G}$-conjugation if and only if 
$A(\bfx) \subset \wt{G}$ is linearly reductive.
By Prop. \ref{prop: finiAAbar2} this happens if and only if $\ov A(\bfx)$ is linearly reductive.
Together with Lemma~\ref{lem: action cmpts} this
implies the assertion.
\epfms

Note that in general it really is necessary to work with  $\ov {A}(\bfx)$ 
rather than ${A}(\bfx)$:

\begin{lem}
There are examples of  $\bfx\in X$
which are polystable but ${A}(\bfx)$
is not linearly reductive.
\end{lem}

\pf 
Write $x_i=(g_i,\phi_i)\in G(\phi_i)$.
Choose $g_1,\ldots,g_n\in G$ generating a Zariski dense subgroup of an abelian
unipotent subgroup $U\subset G$, such as a root group.
Then define $\phi_i$ to be the inner automorphism  
$\phi_i(g)=g^{-1}_i g g_i$.
It follows that ${A}(\bfx)\cong U$ 
is not reductive, but 
$\ov {A}(\bfx)=\{1\}$ is trivial (and so $\bfx$ is polystable by Thm. \ref{thm:R1-tw}).
\epf

\begin{thm}\label{thm:R2-tw}
For any point $\bfx \in X$ the following hold:
\begin{enumerate}
\item $\bfx$ is polystable if and only if 
any ${A}(\bfx)$-invariant parabolic 
subgroup $P \subset G$ 
has an ${A}(\bfx)$-invariant 
Levi subgroup $L \subset P$. 
\item $\bfx$ is stable if and only if 
there are no proper ${A}(\bfx)$-invariant parabolic subgroups of $G$. 
\end{enumerate}
\end{thm}

We prepare three lemmas.

\begin{lem}\label{lmm:levi}
Let $P$ be a parabolic subgroup of $G$ and $L$ be a Levi subgroup of $P$.
Then $N_P(L) = L$.
\end{lem}

\pf
Let $\pi \colon P \to P/R_u(P)$ be the quotient of $P$ by the unipotent radical.
Suppose that $g \in P$ normalises $L$ and decompose it 
as $g=hu$, $h \in L$, $u \in R_u(P)$. 
Then for any $x \in L$ we have $gxg^{-1} \in L$ and hence $uxu^{-1} \in L$.
Since $u \in \Ker \pi$ we have $\pi(uxu^{-1})=\pi(x)$, 
which implies $uxu^{-1}=x$ since the restriction of $\pi$ to $L$ is injective.
Thus $u$ commutes with $L$. 
On the other hand, $L$ coincides with its centraliser 
since $L=C_G(S)$ for some torus $S \subset L$.
Thus $u \in L$ and hence $u=1$, i.e. $g \in L$.
\epf

\begin{lem}\label{lmm:kernel}
The kernel of the $G$-action on $X$ is equal to $Z(G)^\Gamma$.
\end{lem}

\pf
If $k \in G$ lies in the kernel, then for any $g \in G$ and $i=1,2, \dots ,n$ 
we have $k(g,\phi_i)=(g,\phi_i)k$, i.e. 
\[
k g = g \phi_i(k).
\]
Taking $g$ to be $1$ we obtain $\phi_i(k)=k$ ($i=1,2, \dots ,n$). 
Thus $k g = g k$ for all $g \in G$, i.e., $k \in Z(G)$.
Since $\Gamma$ is generated by $\Inn(G), \phi_1,\phi_2, \dots ,\phi_n$  
and $\Inn(G)$ acts trivially on $Z(G)$ we obtain $k \in Z(G)^\Gamma$.
The converse is clear.
\epf

\begin{lem}\label{lmm:ker-stab}
$Z(G) \cap G_\bfx = Z(G)^\Gamma$ for any $\bfx \in X$.
\end{lem}

\pf
Take any $\bfx =(x_i)_{i=1}^n \in X$.
For each $i=1,2, \dots ,n$ we have 
$\Ad(x_i)\phi_i^{-1} \in \Inn(G)$, which acts trivially on $Z(G)$.
Thus for any $g \in Z(G) \cap G_\bfx$ we have
\[
\phi_i(g) = \Ad(x_i)(g) = g \quad (i=1,2, \dots ,n),
\]
which implies $g \in Z(G)^\Gamma$. 
The inclusion $Z(G)^\Gamma \subset Z(G) \cap G_\bfx$ is clear. 
\epf

\pfms (of Theorem~\ref{thm:R2-tw}).\ 
As in the proof of Theorem~\ref{thm:R1-tw}, 
we may assume that $\phi_i$, $i=1,2, \dots ,n$ are 
all contained in a common finite subgroup 
$\Gamma' \subset \Gamma$ 
such that $\Gamma = \Inn(G)\cdot \Gamma'$ and 
put $\wt{G} = G \ltimes \Gamma'$.
Under this assumption $\ov{A}(\bfx)$ is linearly reductive if and only if 
$A(\bfx) \subset \wt{G}$ is linearly reductive,
by Prop.~\ref{prop: finiAAbar2}.

(1) We should show that a subgroup 
$H\subset \wt G$ is linearly reductive if and only if
any $H$-invariant parabolic in $G$ has an 
$H$-invariant Levi subgroup.
First suppose that $H\subset \wt G$ is linearly reductive
and $P \subset G$ is 
an $H$-invariant parabolic subgroup.
Since $H$ is linearly reductive and contained in 
$N_{\wt{G}}(P)$, 
it is contained in some Levi subgroup 
$\wt{L}$ of $N_{\wt{G}}(P)$.
Then $L:=\wt{L} \cap G$ is a Levi subgroup of  
$N_{\wt{G}}(P) \cap G = N_G(P)= P$ 
and normalised by $H$ (so it is $H$-invariant).
Conversely, suppose 
any $H$-invariant parabolic in $G$ has an 
$H$-invariant Levi subgroup.
By \cite[Prop.~2.6]{Richardson}, 
there exists a one-parameter subgroup $\lambda$ of $G$ such that 
$H \subset P_{\wt{G}}(\lambda)$ and  
$R_u (H) \subset U_{\wt{G}}(\lambda)$, 
where 
\begin{align}
P_{\wt{G}}(\lambda) &= \{\, x \in \wt{G} \mid 
\text{$\lim_{t \to 0} \lambda(t) x \lambda(t)^{-1}$ exists}\,\}, \\
U_{\wt{G}}(\lambda) &= \{\, x \in \wt{G} \mid 
\lim_{t \to 0} \lambda(t) x \lambda(t)^{-1} =1 \,\},
\end{align}
and $R_u(H)$ 
is the unipotent radical of $H$.
Put $P = P_G(\lambda) = P_{\wt{G}}(\lambda) \cap G$, 
which is a parabolic subgroup of $G$. 
Since $P$ is normalised by $P_{\wt{G}}(\lambda)$, 
it is normalised by $H$.
Hence $P$ has an $H$-invariant Levi subgroup $L$ by assumption. 
We have $H \subset N_{\wt{G}}(L)$ and hence 
$R_u(H) \subset N_G(L)$ 
(recall that the unipotent radical is connected). 
Note that $R_u(H)$ is also contained in 
$U_{\wt{G}}(\lambda) \cap G = R_u(P)$, 
while any non-trivial element of $R_u(P)$ does not normalise the Levi subgroup $L$
by Lemma~\ref{lmm:levi}. 
Thus $R_u(H)$ is trivial, i.e.  $H$ 
is linearly reductive.

(2) Suppose that $\bfx=(x_i)_{i=1}^n \in X$ is stable and 
let $P$ be a ${A}(\bfx)$-invariant 
parabolic subgroup of $G$.
By Theorem~\ref{thm:R1-tw}, $\ov{A}(\bfx)$ 
(and hence $A(\bfx)$)  
is linearly reductive. 
Since $A(\bfx)$ normalises $P$, 
there exists a Levi subgroup $\wt{L}$ of $N_{\wt{G}}(P)$ containing $A(\bfx)$.
By \cite[Prop.~2.4]{Richardson}, there exists a one-parameter 
subgroup $\lambda$ of $G$ such that  
$N_{\wt{G}}(P) = P_{\wt{G}}(\lambda)$ and  
$\wt{L}=C_{\wt{G}}(\range \lambda)$.
Since $A(\bfx) \subset C_{\wt{G}}(\range \lambda)$ each $x_i$ 
commutes with $\range \lambda$ and hence 
$\range \lambda \subset G_\bfx$.
The stability now implies that $\range \lambda$ is contained in the kernel 
$Z(G)^\Gamma$ and hence $P=G$.
Conversely, suppose that $\bfx \in X$ is not stable.
Then if the orbit $G \cdot \bfx$ is not closed 
we argue as follows:
By the Hilbert--Mumford criterion, 
there exists a one-parameter subgroup $\lambda$ of $G$ and 
an element $\bfy=(y_i)_{i=1}^n \in X$
such that $\lim_{t \to 0} \lambda(t) \cdot \bfx = \bfy$ and 
$G \cdot \bfy$ is closed.
We show that the parabolic subgroup 
$P := P_G(\lambda) = P_{\wt{G}}(\lambda) \cap G$
is ${A}(\bfx)$-invariant and proper.
For any $g \in P$ and $i=1,2, \dots ,n$, the limit of 
\[
\lambda(t) \Ad(x_i)(g) \lambda(t)^{-1} 
= \lambda(t) x_i \lambda(t)^{-1} \cdot 
\lambda(t)g \lambda(t)^{-1} \cdot 
\lambda(t) x_i^{-1} \lambda(t)^{-1}
\]
as $t \to 0$ exists. 
Hence $P$ is ${A}(\bfx)$-invariant.
If $P$ is not proper, 
$\range \lambda$ is contained in 
$Z(G) \cap G_\bfy =Z(G)^\Gamma$.
Thus $\lambda(t) \cdot \bfx = \bfx$ ($t \in \IC^*$) and hence 
$\bfx = \bfy$, which contradicts 
the assumption that $G \cdot \bfx$ is not closed.
Hence $P$ is proper. 
Finally suppose $G \cdot \bfx$ is closed, but of the wrong dimension.
Then the stabiliser $G_\bfx$ is linearly reductive 
(\cite{Richardson} 1.3.3) and 
the quotient $G_\bfx/Z(G)^\Gamma$ has non-trivial
identity component.
Hence there exists a one-parameter subgroup $\lambda$ of $G_\bfx$ 
such that $\range \lambda \not\subset Z(G)^\Gamma$.
Since each $x_i$ commutes with $\lambda$,  
the parabolic subgroup $P:=P_G(\lambda)$ is ${A}(\bfx)$-invariant. 
It is proper since Lemma~\ref{lmm:ker-stab} implies 
$\range \lambda \not\subset Z(G)$.
\epfms

Let us rephrase/abstract the first part of Theorem~\ref{thm:R2-tw} in a way that will be useful later.
Let $G$ be a reductive group and $\Lambda$ an algebraic group 
acting on $G$ by (algebraic) group automorphisms.
Suppose that the action of $\Lambda$ is effective and 
the identity component $\Lambda^0$ acts by inner automorphisms.
Then we may regard $\Lambda$ as a subgroup of $\Aut(G)$ 
and its image in $\Out(G)$ is finite.

\begin{prop}\label{prop: abstr-compred}
$\Lambda$ is linearly reductive if and only if 
any $\Lambda$-invariant parabolic subgroup of $G$ has a 
$\Lambda$-invariant Levi subgroup.
\end{prop}

\pf
By the result of Borel--Serre and Brion, 
there exists a finite subgroup $\Lambda' \subset \Lambda$ such that 
$\Lambda^0 \Lambda' = \Lambda$.
Put $\wt{G} = G \ltimes \Lambda'$ and let 
$K \subset \wt{G}$ 
be the preimage of $\Lambda$ under $\Ad$,
 so that Prop.~\ref{prop: finiAAbar2} implies 
$\Lambda$ is linearly reductive if and only if $K$ 
is linearly reductive.
Thus the equivalence follows from 
the proof of Thm.~\ref{thm:R2-tw} (1).
\epf

In the same setting there is a further characterisation:

\begin{prop} \label{prop: Lred-levi}
$\Lambda$ is linearly reductive if and only if 
there exists a torus $S \subset G$ such that $C_G(S)$ 
is $\Lambda$-invariant 
and has no proper $\Lambda$-invariant parabolic subgroups.
\end{prop}

\pf
Suppose that $\Lambda$ is linearly reductive. 
We show that there exists a decreasing sequence 
of $\Lambda$-invariant closed subgroups
\[
G=L_0 \supset P_1 \supset L_1 \supset P_2 \supset L_2 \supset 
\cdots \supset P_r \supset L_r
\]
such that each $P_i$ is a proper parabolic subgroup of $L_{i-1}$, 
each $L_i$ ($i>0$) is a Levi subgroup of $P_i$ and 
$L_r$ has no proper $\Lambda$-invariant parabolic subgroups. 
If $L_0=G$ has no proper $\Lambda$-invariant parabolic subgroups 
we have nothing to do (just put $r=0$).
Otherwise we take any proper $\Lambda$-invariant parabolic subgroup 
$P_1 \subsetneq L_0$. Then by Prop. \ref{prop: abstr-compred} 
there exists a $\Lambda$-invariant Levi subgroup $L_1 \subset P_1$.
Let $\Lambda_1$ be the quotient of $\Lambda$ by the kernel of 
the induced $\Lambda$-action on $L_1$, so that 
$\Lambda_1$ effectively acts on $L_1$. 
Note that its identity component $\Lambda_1^0$  
acts by inner automorphisms of $L_1$; 
indeed, the action of any element of $\Lambda_1^0$ is induced from 
some inner automorphism of $G$ preserving $P_1, L_1$  
and hence is inner by Lemma~\ref{lmm:levi}. 
Since $\Lambda$ is linearly reductive $\Lambda_1$ is also linearly reductive, 
and a subgroup of $L_1$ is $\Lambda_1$-invariant 
if and only if it is $\Lambda$-invariant as a subgroup of $L_0$.
If $L_1$ has no proper $\Lambda_1$-invariant parabolic subgroup, 
the sequence $L_0 \supset P_1 \supset L_1$ is as desired.
Otherwise we take any proper $\Lambda_1$-invariant parabolic subgroup 
$P_2 \subsetneq L_1$. Then by Prop. \ref{prop: abstr-compred}  
there exists a $\Lambda_1$-invariant Levi subgroup $L_2 \subset P_2$.
Iterating this procedure, we obtain a desired decreasing sequence.
Since each $L_i$ ($i>0$) is a Levi subgroup of $P_i$ 
there exists a torus $S_i \subset L_i$ such that $L_i = C_{L_{i-1}}(S_i)$.
Then $L_r$ is the common centraliser of $S_1, S_2, \dots , S_r$ in $G$.
Hence the torus $S \subset G$ generated by $S_1, S_2, \dots , S_r$ (note that they 
commute with each other) is as desired.

Conversely, suppose that there exists a torus $S \subset G$ 
such that $L:=C_G(S)$ is $\Lambda$-invariant 
and has no proper $\Lambda$-invariant parabolic subgroups. 
Let $P \subset G$ be a $\Lambda$-invariant parabolic subgroup. 
Then the intersection
$P \cap L$ is a $\Lambda$-invariant parabolic subgroup of $L$ 
and hence $P \cap L = L$, i.e. $L \subset P$.
Since $L$ is reductive, there exists a Levi subgroup $M \subset P$ containing $L$. 
For any $\psi \in \Lambda$ the image $\psi(M)$ of $M$ is 
also a Levi subgroup of $P$ containing $L$. Since $L$ 
contains a maximal torus of $P$ and any maximal torus of $P$ is contained 
in a unique Levi subgroup, we have $\psi(M) = M$. 
Hence $M$ is $\Lambda$-invariant, which together with Prop. 
\ref{prop: abstr-compred} 
shows that $\Lambda$ is linearly reductive.
\epf

Applying this to the set-up  of the present section
(with $\La=\ov A(\bfx)$) yields:

\begin{cor}\label{cor: completered}
The following are equivalent:

0) A point $\bfx \in X$ is polystable,

1) The group $\ov A(\bx)$ is linearly reductive,

2) Any ${A}(\bfx)$-invariant parabolic in $G$
has an ${A}(\bfx)$-invariant Levi subgroup, 

3) There exists a subtorus $S \subset G$ such that 
$C_G(S)$ is ${A}(\bfx)$-invariant and has no proper 
${A}(\bfx)$-invariant parabolic subgroups,

4) 
There exists an ${A}(\bfx)$-invariant Levi subgroup 
$L$ of a parabolic of $G$, such that $L$ has no proper 
${A}(\bfx)$-invariant parabolic subgroups.

\end{cor}

Note that 3) and 4) are trivially equivalent since 
centralisers of tori in $G$ are exactly 
the Levi subgroups of parabolics.

\section{More general set-up}\label{sn: moregen}

The results of the previous section will now be generalised, in a form more directly useful in the context of Stokes local systems.  
Return to the set-up of Thm. \ref{thm:R1-tw} (with $n\ge 1$), but now 
further choose an integer $m\ge 1$ and tori
$\IT_1,\ldots,\IT_m\subset G$.
Write $\bT=\IT_1\times\cdots\times \IT_m\subset G^m$,
let $H_i=C_G(\IT_i)$ and 
$\bH=H_1\times\cdots \times H_m
=C_{G^m}(\bT)\subset G^m$.
We allow some of the $\IT_i$ to be a point, in which case $H_i=G$.
In this section we will study the stability and polystability for the action of $\bH$ on
\beq\label{eq: twRich}
X:=G^{m-1}\times \Prod_{i=1}^nG(\phi_i)
\eeq
given by 
$$\bh\cdot(\bC,\bM) = 
(h_2C_2h_1^{-1},\ldots,h_mC_mh_1^{-1},
h_1M_1h_1^{-1},\ldots,h_1M_nh_1^{-1})
$$
where $\bh = (h_1,\ldots,h_m),\bC=(C_2,\ldots,C_m), \bM=(M_1,\ldots,M_n), C_i\in G, M_i\in 
G(\phi_i).$

Thus if $m=1$ and $\IT_1=1$ we recover the situation
of Thm. \ref{thm:R1-tw}.
The case $m=1$ and $\IT_1$ arbitrary but each $\phi_i=1$ was studied by Richardson in \cite{Richardson} Thm. 13.2,14.1
(taking $S=\IT_1$ acting on $G$ by conjugation).
More generally, in effect,
\cite{gbs} Cor. 9.6 studied
the notion of stability in 
the case with $m$ arbitrary and each $\phi_i=1$,
making the link to the differential Galois group 
of irregular connections (whence the $\IT_i$ are the Ramis/exponential tori).

For $\bx=(\bC,\bM)\in X$ let 
$A(\bx)\subset G\ltimes \Ga$ be the 
Zariski closure of the subgroup generated by
$$M_1,M_2,\ldots,M_n,\IT_1,C_2^{-1}\IT_2C_2,\ldots,
C_m^{-1}\IT_mC_m.$$
Let $\ov A(\bx)\subset \Ga$ be the image of $A(\bx)$ in $\Ga\subset \Aut(G)$, as before.
Recall that a subset of $G$ is $A(\bx)$-invariant if it is preserved by $\ov A(\bx)\subset \Aut(G)$.

\begin{thm}\label{thm: moregen}
(1) A point $\bx\in X$ is polystable for the $\bH$ action if and only if $\ov A(\bx)$ is linearly reductive.

(2) A point $\bx\in X$ is stable for the $\bH$ action if and only if there are no proper $A(\bx)$-invariant
parabolic subgroups in $G$.
\end{thm}

As in the last section, one can rephrase polystability
in several different ways 
(recall Propositions \ref{prop: abstr-compred}, \ref{prop: Lred-levi}):

\begin{cor}\label{cor: completered2}
A point $\bfx \in X$ is polystable if and only if

1) The group $\ov A(\bx)$ is linearly reductive, or

2) Any ${A}(\bfx)$-invariant parabolic in $G$
has an ${A}(\bfx)$-invariant Levi subgroup, or

3) There exists a subtorus $S \subset G$ such that 
$C_G(S)$ is ${A}(\bfx)$-invariant and has no proper 
${A}(\bfx)$-invariant parabolic subgroups.
\end{cor}

\pfms (of Thm. \ref{thm: moregen}).
Let $G\times \bH$ act on  $\wt X:= G^m\times \Prod_1^N G(\phi_i)$
via 
$$(g,\bh)\cdot(\bC,\bM) = 
(h_1C_1g^{-1},h_2C_2g^{-1},\ldots,h_mC_mg^{-1},
gM_1g^{-1},\ldots,gM_ng^{-1})
$$
where $\bh = (h_1,\ldots,h_m),\bC=(C_1,\ldots,C_m), \bM=(M_1,\ldots,M_n), C_i\in G, M_i\in G(\phi_i).$
(Up to relabelling and incrementing $m$
this is the special case where $H_1=G$.)
Consider the $\bH$-equivariant map $\wt X \to X$ taking $(\bC;\bM)$ to
$$(C_2C_1^{-1},\ldots, C_mC_1^{-1};
C_1M_1C_1^{-1},\ldots,C_1M_NC_1^{-1})\in X.$$
It expresses $\wt X$ as a principal $G$-bundle over $X$ (the fibres are exactly the $G$-orbits).
Any point $\bx\in X$ has a unique lift $\wt \bfx\in \wt X$ with $C_1=1$. 
It follows that the $\bH$ orbit of $\bx\in X$ is closed if and only if the $G\times \bH$ orbit of $\wt\bfx$ is closed in $\wt X$.

Now choose $\bt=(t_1,\ldots,t_m)\in \bT$ so that $t_i$ generates a Zariski dense subgroup of $\IT_i$ for each $i$.
In particular $H_i=C_G(t_i)$.
Consider the simultaneous conjugation action of 
$G$ on $Y:=G^m\times \Prod_1^N G(\phi_i)$, and the $G$-equivariant embedding
$$\pi:\wt X/\bH\into Y;\quad [(\bC;\bM)]\mapsto (\bC^{-1}\bt \bC,\bM).$$
Thus $\wt X/\bH$ is identified with a closed subvariety of 
$Y$ (by \cite{Richardson} p.1),  
and $\wt X$ is a $G$-equivariant principal $\bH$-bundle over 
$\wt X/\bH$.
It follows that the
$G\times \bH$ orbit of $\wt\bfx$ is closed in $\wt X$
if and only if the $G$ orbit of $\pi(\wt \bx)$ is closed in $Y$.

Hence part (1) of the theorem follows from Thm. \ref{thm:R1-tw} 
(applied to $Y$).
To deal with stability we need to consider the stabilisers and the kernels of the actions.

\begin{lem}\label{lem: stabcorr}
The stabiliser of any point $\bx \in X$ for the $\bH$-action
is canonically isomorphic to the stabiliser of the point 
$\pi(\wt \bx)\in Y$ for the $G$-action.  
\end{lem}

\pf
Suppose $g \in G$ fixes $(\bC^{-1} \bt \bC, \bM)$
and let
$h_i = C_i g C_i^{-1}$ ($i=1,2, \dots ,m$) and 
$\bh =(h_1,\ldots,h_m)$.
Note  $h_1=g$ because $C_1=1$.
Since $g$ commutes with each $C_i^{-1} t_i C_i$, 
each $h_i$ commutes with $t_i$, which implies $\bh \in \bH$.
We have $h_i C_i g^{-1} = C_i$ ($i=1,2, \ldots ,m$) by the definition of $h_i$ and 
$g M_j g^{-1} = M_j$ ($j=1,2, \dots ,N$) as $g$ centralises $\bM$.
Hence the pair $(g,\bh) \in G \times \bH$ stabilises 
$\wt \bx \in \wt X$, 
and hence $\bh$ stabilises $\bx \in X$.
Conversely if $\bh =(h_1,\ldots,h_m) \in \bH$ 
fixes $\bx$, then let $g = h_1$. 
It follows immediately that  $(g,\bh)\in G\times \bH$ 
fixes $\wt \bx$ and thus that $g$ fixes 
$\pi(\wt \bx)$.
Clearly the two correspondences are inverses of each other.
\epf

As in Lemma \ref{lmm:kernel} one has:
\begin{lem}
{\rm (1)} The kernel of the $G$-action on $Y$ (or $\wt X/\bH$)
is the $\Gamma$-invariant subgroup $Z(G)^\Gamma$ of the center of $G$, and
{\rm (2)} The kernel of the $\bH$-action on $X$ 
is the subgroup
of elements $(h_1,\ldots,h_m) \in \bH$ satisfying
$h_1 = h_2 = \cdots = h_m \in Z(G)^\Gamma$.
\end{lem}

Note that these two groups correspond to each other under the correspondence of Lem. \ref{lem: stabcorr} (for any $\bx\in X$).
It follows that $\bx\in X$ is stable if and only if $\pi(\wt \bx)\in Y$ is stable. By Thm. \ref{thm:R2-tw} (2) applied to $Y$, 
this happens if and only if there are no proper $A(\bx)$ invariant parabolic subgroups of $G$.
\epfms

\section{Application to wild character varieties}\label{sn: apwcv}

Recall from \cite{gbs,twcv} that 
the wild character variety 
$\MB(\bSi)=\THom_\IS(\Pi,G)/\bH$ 
is determined by an 
irregular curve/wild Riemann surface 
$\bSi=(\Si,\al,\Th)$ with group $\cG$,  
where $\Si$ is a compact smooth complex algebraic curve, 
$\al\subset \Si$ is a non-empty finite subset, 
$\cG\to \Si^\circ:=\Si\setminus \al$ is a local system of groups over the punctured curve (with each fibre isomorphic to some fixed connected complex reductive group $G$), 
and $\Th$ consists of the data of an irregular class 
$\Th_a$ at each point $a\in\al$ 
(in the sense of \cite{twcv} \S3.5---it is the class of a graded $\cG$ local system).

As in \cite{twcv} \S4, 
$\bSi$ then determines an auxiliary surface 
$\wt \Si$, equipped with boundary circles $\partial$,
halos $\IH\subset\wt \Si$,
and tangential punctures $e(\IA)$ 
(near the  singular directions $\IA$). 
Further $\bSi$ determines 
a local system $\IT\to \partial$ of finite dimensional complex 
tori, the Ramis tori (\cite{twcv} p.9).
Choosing a finite set 
$\be=\{b_1,\ldots,b_m\}
\subset \partial$ of basepoints 
(with one point in each component circle, 
as in \cite{twcv} \S4.1) then determines 
the wild surface groupoid  $\Pi=\Pi_1(\wt \Si,\be)$,
the fundamental groupoid of the auxiliary surface with 
these basepoints, as in \cite{twcv} \S5. 
The local system $\cG$ is determined by a map 
$f:\Pi\to \Aut(G)$ and 
this determines the space of $f$-twisted representations
$$\THom(\Pi,G) =\{\rho\in \Hom(\Pi,G\ltimes \Aut(G))\st 
\rho(\ga)\in G(f(\ga))\text{ for all }\ga\in \Pi\}$$ 
of $\Pi$, as in \cite{twcv} \S5.
Since $\wt \Si$ is just a punctured real surface with boundary, 
choosing generating paths in $\Pi$ yields 
an isomorphism $\THom(\Pi,G)\cong G^N$ of spaces, for some integer $N$, so it is a smooth affine variety.
In turn the wild representation variety 
$\cR=\THom_\IS(\Pi,G)$ is the closed subvariety of
$\THom(\Pi,G)$ 
cut out by the two Stokes conditions, in \cite{twcv} Defn. 18.
Intrinsically, $\cR$ is the moduli space of framed Stokes local systems, as in \cite{twcv} Prop. 19, framed via a graded isomorphism to a standard fibre $\cF_i$ at each basepoint $b_i\in \be$ 
(\cite{twcv} \S4.1).
The group $H_i=\GrAut(\cF_i)=C_G(\IT_i)\subset G=\Aut(\cF_i)$ acts transitively on  the set of framings at $b_i$, where 
$\IT_i=\IT_{b_i}$. 
Thus the group  $\bH=\Prod_{b_i\in \be} H_i$ acts naturally on $\cR$.
The wild character variety $\MB(\bSi)$ is the affine 
geometric invariant theory quotient $\cR/\bH$, and so its points are the closed $\bH$ orbits in $\cR$.
This leads directly to the key statement:

\begin{prop}\label{prop: embed-wcv-pres}
The wild representation variety $\cR=\THom_\IS(\Pi,G)$
may be embedded in an $\bH$-equivariant way, as a closed
subvariety of the space $X=G^{m-1}\times \Prod_1^n G(\phi_i)$ 
of \eqref{sn: moregen}, for suitable $n$ and automorphisms 
$\{\phi_i\}\subset\Aut(G)$, with $m=\#\al$.
\end{prop}
\pf
This comes down to considering the inclusion
$\THom_\IS(\Pi,G)\subset \THom(\Pi,G)$ as a closed subvariety 
(forgetting the Stokes conditions, as in \cite{twcv} Defn. 18), 
and then identifying $\THom(\Pi,G)\cong X$
in an $\bH$-equivariant way.
Both of these are straightforward.
In particular 
$M_j=\rho(\ga_j)$ for generators $\ga_j$
of $\pi_1(\wt\Si,b_1)$, and $C_i=\rho(\chi_i)$ for paths $\chi_i$ from $b_1$ to $b_i$, for $i=2,\ldots,m$. 
\epf

\begin{rmk}
Although we do not directly need it here, 
note that the Stokes $G$-local 
systems on $\bSi$ encode the algebraic connections on algebraic principal  $G$-bundles on the open curve 
$\Si^\circ:= \Si\setminus \ba$, with 
irregular class $\Th_a$ at each $a\in \ba$.
This follows from the Tannakian approach of \cite{MR91}
once you identify $\rho:\Pi\to G$ (plus the irregular class) 
as giving a  
representation of the wild fundamental group of \cite{MR91}.
  \footnote{\cite{MR91} is somewhat sketchy, due to the untimely passing of one of the authors. 
However we do indeed live in the sunlit uplands described there, since, if we use (1), (2), (3) to label the approaches of Malgrange--Sibuya, Deligne--Malgrange, and Martinet--Ramis
respectively,
then, in effect,  
\cite{Mal-irregENSlongplusnote} proves (1) $\Leftrightarrow$ (2), \cite{L-R94} proves (1) $\Leftrightarrow$ (3), and 
\cite{tops} gives a direct proof (2) $\Leftrightarrow$ (3) avoiding (1), passing via (4) the Stokes graded local systems.}
\end{rmk}

\begin{rmk}
We don't worry about 
the exact positions of the tangential punctures
because one always gets the same groupoid $\Pi$;
as in \cite{gbs} we attach the tangential punctures
to the singular directions $\IA\subset \d$ via noncrossing cilia
and then shrink the cilia whenever convenient, identifying $\Pi$ 
in the limit with that defined from the canonical, 
non-Hausdorff, version of $\wt \Si$ as in \cite{tops} Rmk 8.4, 2).
\end{rmk}

\subsection{Galois groups from Stokes representations}\label{ssn: gg from sreps}
In order to 
define the Galois group $\Gal(\rho)$ of $\rho\in \THom_\IS(\Pi,G)$
we will assume that $\cG$ is ``Out-finite'', in the sense 
that the monodromy group 
$f(\pi_1(\wt\Si,b_1))\subset \Aut(G)$ of $\cG$ 
has finite image in $\Out(G)=\Aut(G)/\Inn(G)$.
Then the group $\Ga$ generated by $f(\pi_1(\wt\Si,b_1))$ and $\Inn(G)$ 
is an algebraic group, as in \S\ref{sn: trich} above. 
Thus any Stokes representation takes values in the algebraic group 
$G\ltimes \Ga\subset G\ltimes \Aut(G)$,
and so we can consider the Zariski closure of its monodromy.
The Galois group is defined by adding the Ramis tori as well:
If $t\in \IT_i\subset G$ and 
$\chi$ is any path in $\wt \Si$ from $b_1$ to $b_i$, 
consider the element 
\beq\label{eq: toriinGv1}
C^{-1}tC\in G=G(\Id)\subset G\ltimes \Ga,
\eeq 
where $C=\rho(\chi)$.

\begin{defn}\label{defn: galrho}
The differential Galois group $\Gal(\rho)$ of $\rho$
is the Zariski closure of the subgroup of $G\ltimes \Ga$
generated by $\rho(\pi_1(\wt\Si,b_1))$ and all of the tori 
\eqref{eq: toriinGv1} (as $i,t,\chi$ vary).
\end{defn}

It follows that $\Gal(\rho)$ acts on $G$ by group automorphisms, via the adjoint action of $G\ltimes \Ga$ on $G=G(\Id)$.
Let $\ov \Gal(\rho)\subset \Ga\subset \Aut(G)$ be the resulting 
image of $\Gal(\rho)$. 
This definition is (of course) motivated by 
Ramis' description (\cite{ramis-pv}, 
 \cite{MR91} Thm. 21, \cite{L-R94} Thm. III.3.11) 
 of the differential Galois group 
of an algebraic connection on a vector bundle.

\begin{cor}\label{cor: main-red-annonce}
A Stokes representation $\rho\in \THom_\IS(\Pi,G)$ is polystable for the action of $\bH$ if and only if $\bar\Gal(\rho)$ is a linearly reductive group.
\end{cor}
\pf
This now follows from part (1) of Thm. \ref{thm: moregen}, via Prop. 
\ref{prop: embed-wcv-pres}, since  $\bar\Gal(\rho)$ matches up with 
$\ov A(\bx)$. 
\epf

Special cases include:

$\bullet$ If $\cG$ has finite monodromy then
$\rho\in \THom_\IS(\Pi,G)$ is polystable if and only if $\Gal(\rho)$ is a linearly reductive group.

$\bullet$ If $\cG$ is a constant general linear group then
$\rho\in \Hom_\IS(\Pi,G)$ is polystable if and only if
$\rho$ is the direct sum of irreducible Stokes representations.

Recall that $Z(G)^\Ga$ is the $\Ga$ invariant subgroup of the centre of $G$, and it embeds diagonally in 
$\bH$.
To deal with stability 
we will avoid degenerate cases
by assuming:
\beq \label{eq: stokes gluing}
\begin{matrix}
\text{\em The kernel of the action of $\bH$ on 
$\THom_\IS(\Pi,G)$  is $Z(G)^\Ga$. }
\end{matrix}
\eeq
The lemma below shows one can always add one or two punctures to ensure this condition holds.
Note that no generality is lost: any symplectic leaf 
$\MB(\bSi,\cC)\subset \MB(\bSi)$ will also be a symplectic leaf of the larger Poisson variety obtained by first making such additional punctures (namely that with trivial monodromy around the new punctures).
For example 
the usual character variety of a genus $g>0$ compact Riemann surface $\Si$  is a (very special) symplectic leaf of the character variety of $\Si\setminus a$, for any point $a\in \Si$.

\begin{lem} \label{lem: suff kernel condition}
Suppose $m\ge 1$ and $a_1\in \al$ has trivial irregular class, 
and if $g=0$ then $m\ge 2$.
Then $\THom_\IS(\Pi,G)$ is a smooth non-empty affine variety and the kernel $\bK$ of the $\bH$ action is $Z(G)^\Ga$ embedded 
diagonally in $\bH$.
\end{lem}
\pf
It is nonempty as it is the fusion of some fission spaces and some internally fused doubles:
Recall that $\THom_\IS(\Pi,G)$ can be described as 
the quasi-Hamiltonian $G$-reduction 
\beq
\THom_\IS(\Pi,G) \cong 
\bigl(\ID_1\fusion{G}\cdots \fusion{G} \ID_g \fusion{G}\cA(Q_1)\fusion{G}\cdots \fusion{G}\cA(Q_m)\bigr)\spq G
\eeq
where each $\ID_i$ is a (twisted) internally fused double (\cite{twcv} p.23, \cite{gbs} Thm 8.2).
Since $\cA(Q_1)\cong \bD(G)$ is the double of $G$, it follows that 
$$\THom_\IS(\Pi,G) \cong 
\ID_1\fusion{}\cdots \fusion{} \ID_g \fusion{}\cA(Q_2)\fusion{}\cdots \fusion{}\cA(Q_m)$$ 
so that $\THom_\IS(\Pi,G)$  is the product of some smooth nonempty affine varieties.
Note that $\bH$ still acts on $\THom_\IS(\Pi,G)$ 
and this includes $H_1=G$.

Suppose $\bk\in \bK$ and suppose (as usual) that the framings of $\cG$ are such that the monodromy in $\Aut(G)$ of $\cG$ is trivial along the $m-1$ 
chosen paths $\chi_i:b_1\to b_i$. 
Then $\bk$ acts on $C_i=M_{\chi_i}$ as $k_iC_ik_1^{-1}$.
Taking $C_2=C_3=\cdots C_m=1$ implies $k_1=k_2=\cdots =k_m$.
If $m\ge 2$ then the fact that all $C_2\in G$ are fixed implies
$k_1\in Z(G)$.
On the other hand if $m=1$ then $g>0$ so 
looking at $\ID_1$ 
there is $\phi\in \Aut(G)$ so that 
$k_1 A \phi(k_1^{-1})=A$ for all $A\in G$.
This implies $k_1 =\phi(k_1)$ and $k_1\in Z(G)$.
Thus in all case $k_1$ is central. 
Then looking at any loop $\ga$ based at $b_1$ leads to a relation of the form $k_1 M_\ga \phi_\ga(k_1^{-1})=M_\ga$.
Thus since $k_1$ is central this implies 
$k_1 = \phi_\ga(k_1)$, so $k_1\in Z(G)^\Ga$.
\epf

\begin{rmk}
Note that it follows in general  (as in \cite{gbs} Thm. 8.2) that if $\THom_\IS(\Pi,G)$
is nonempty (and $m>0$) then it is a smooth
affine variety.
\end{rmk}

Part (2) of Theorem \ref{thm:R2-tw} then implies:

\begin{cor}\label{cor: irredGalrho}
Suppose the condition \eqref{eq: stokes gluing} holds 
(for example if we make additional punctures so the conditions 
in Lemma \ref{lem: suff kernel condition} hold).
Then a Stokes representation $\rho\in \THom_\IS(\Pi,G)$ is stable for the action of 
$\bH$ if and only if there is no proper parabolic subgroup $P\subset G$ stabilised by 
the action of $\Gal(\rho)$.
\end{cor}
\pf
This follows from Thm. \ref{thm:R2-tw}
since $\cR$ is closed in $X$, 
$\Gal(\rho)$ matches up with $A(\bx)$,
and the kernel of the $\bH$ action on $\cR$ and on $X$ is the same.
\epf

\section{Stability and polystability of Stokes local systems}\label{sn: SLS}

This section will consider the intrinsic objects (Stokes local systems) underlying Stokes representations, and define the notions of ``irreducible'' and ``reductive'' for Stokes local systems, in terms of reductions of structure group.
Then we will deduce:

\begin{thm}\label{thm: red-irred-SLS}
Suppose $\IL$ is a Stokes local system and
$\rho\in \THom_\IS(\Pi,G)$ is the monodromy of $\IL$.

1) $\rho$ is stable for the action of $\bH$ if and only if $\IL$ is irreducible,

2) $\rho$ is polystable for the action of $\bH$ if and only if $\IL$ is reductive.
\end{thm}

\subsection{Graded local systems}\label{ss: grlocsys}

A Stokes local system is a special type of $\IT$-graded local system 
(in the sense of Defn. \ref{defn: T-graded LS} below), 
so to clarify the ideas we will focus on them here---the results for Stokes local systems follow almost immediately.

Let $S$ be a connected real oriented surface of finite topological type.
Let $\IH\subset S$ be an open subset, and let $\IT\to \IH$ be a local system of complex tori over $\IH$. We allow the fibres of $\IT$ to have different dimensions in different components of $\IH$.
Fix a connected complex reductive group $G$, and a local system $\cG\to S$ of groups, such that each fibre of $\cG$ is isomorphic to $G$.
We will  assume throughout that $\cG$ is ``Out-finite''. This means that the monodromy of $\cG$ has finite image in $\Out(G)$.
In more detail, given a basepoint $b\in S$ and a framing $G\cong \cG_b$ of $\cG$ at $b$, 
then the monodromy representation 
$f:\pi_1(S,b)\to \Aut(G)$ of $\cG$
is such that the monodromy group $f(\pi_1(S,b))\subset \Aut(G)$ has finite image in $\Out(G)$.
Of course if $G$ is semisimple then $\Out(G)$ is finite and so this is no restriction.

Recall that a $\cG$-local system over $S$ is a local system 
$\IL\to S$ which is a $\cG$-torsor (cf. e.g. \cite{twcv} \S2.1), 
and it determines a local system $\Aut(\IL)\to S$ of groups
(each fibre of which is also isomorphic to $G$).

\begin{defn}\label{defn: T-graded LS}
A $\IT$-graded $\cG$-local system over $S$ is a $\cG$-local system 
$\IL\to S$ together with an embedding 
$$\IT\into \Aut(\IL)\big\vert_\IH$$
of local systems of groups over $\IH$.
\end{defn}

For brevity we will simply call this a ``graded local system on $S$'', 
and write $\IT\into \Aut(\IL)$ for the grading. 
A Stokes local system $\IL$ (in the sense of \cite{gbs, twcv}) 
is a special type of $\IT$-graded local system, taking 
$S=\wt \Si$ to be the auxiliary surface, 
$\IH$ to be the union of the halos, 
and $\IT$ to be the image in $\Aut(\IL)\bigl\vert_\IH$ of the exponential torus $\cT$.
(Note that $\IT$ is determined just by the 
irregular class of $\IL$, in the sense of \cite{twcv} \S3.5, since, as explained there, the class determines the finite rank local system 
$I\subset \cI$ of lattices, and $\IT$ is the local system of tori with character lattice $I$.)

\subsection{Galois group}\label{ss: galois group}

If $\IL\to S$ is a $\IT$-graded local system, 
and $b\in S$ is a basepoint, 
define $\Gal(\IL)$ to be the Zariski closure of the group generated by the monodromy of $\IL$ and all the tori $\IT$
(after transporting them to $b$). 
In more detail, first identify $\cG_b\cong G, \IL_b\cong \cF$ (the trivial $G$-torsor, as in \cite{twcv} \S2) and define 
$\Ga\subset \Aut(G)$ to be the group generated by the monodromy 
of $\cG$ and $\Inn(G)$ (as in \S\ref{sn: trich} above).
If $\ga\in \pi_1(S,b)$ let 
$\rho(\ga)\in G(f(\ga))\subset G\ltimes \Ga\into \Perm(\cF)$
be the monodromy of $\IL$ around $\ga$ (where $\Perm(\cF)$
is the  group of all permutations of the fibre $\cF$).
Similarly if $p\in \IH, t\in \IT_p\subset \Aut(\IL)_p\cong G$ 
and $\chi$ is any path in $S$ from $b$ to $p$, consider the element 
\beq\label{eq: toriinG}
C^{-1}tC\in \Aut(\cF)=G=G(\Id)\subset G\ltimes \Ga,
\eeq 
where $C$ is the transport of $\IL$ along $\chi$. %

\begin{defn}
The differential Galois group $\Gal(\IL)$ of $\IL$
is the Zariski closure of the subgroup of $G\ltimes \Ga$
generated by $\rho(\pi_1(S,b))$ and all of the tori 
\eqref{eq: toriinG} (as $p,t,\chi$ vary).
\end{defn}

It follows that $\Gal(\IL)$ acts on $G$ by group automorphisms, via the adjoint action of $G\ltimes \Ga$ on $G=G(\Id)$.
Let $\ov \Gal(\IL)\subset \Ga\subset \Aut(G)$ be the resulting 
image of $\Gal(\IL)$.
Up to isomorphism the affine algebraic group 
$\Gal(\IL)$ and its action on $G$
do not depend on the choice of basepoint $b$ or framings. 

\subsection{Irreducible graded local systems}\label{ss irred gr loc sys}

Define a graded local system $\IL\to S$ to be {\em reducible}
if $\Aut(\IL)$ has a sublocal system of proper parabolic subgroups, containing $\IT$.
In other words
there is a sublocal system $\cP\subset \Aut(\IL)$ such that
1) each fibre $\cP_b$ is a proper parabolic subgroup of $\Aut(\IL)_b$, and 2) the grading $\IT\into \Aut(\IL)$ factors through $\cP$.
Such $\IL$ is {\em irreducible} if it is not reducible.
\begin{lem}\label{lem: 2irred-sls-defns}
$\IL$ is reducible if and only if $\Gal(\IL)$ preserves a proper parabolic subgroup of $G$ (recalling that $\Gal(\IL)$ naturally acts on $G$ by group automorphisms).
\end{lem}
\pf
Suppose $P\subset G\cong \Aut(\IL_b)$ is preserved by $\Gal(\IL)$.
Then $P$ is the fibre at $b$ of a local system of parabolic subgroups $\cP\subset \Aut(\IL)$, since the monodromy of $\Aut(\IL)$  is given by the adjoint action of the monodromy of $\IL$, which is in 
$\Gal(\IL)$. 
Moreover $\IT\into \cP$, since the transport to 
$b$ of each fibre of $\IT$ is in $G$ and preserves $P$ (so is in $P$, since $N_G(P)=P$).
The converse is similar, taking the fibre at $b$ of 
$\cP\subset\Aut(\IL)$.
\epf

This can be related to (twisted) 
reductions of structure group as follows 
(compare \cite{twcv} Defn. 11).

\begin{defn} Suppose $\IL$ is a $\IT$-graded $\cG$-local system.

\label{defn: twisted reductions of structure}
$\bullet$\ A {\em reduction} of $\IL$ is 
a $\IT$-graded $\cP$-local system
$\IP\to S$  (for some local system of groups $\cP\to S$), such that
$\IL$ is a twisted pushout of $\IP$.
This means that there is a $\cG$ local system $\IM$ with an embedding 
$\cP\into \Aut(\IM)$, together with an isomorphism 
$\IL\cong \IP\times_\cP \IM$ (of graded $\cG$ local systems),

$\bullet$\ A reduction $\IP$ of $\IL$ is a {\em parabolic reduction} if the fibres of $\cP$ embed as parabolic subgroups of the fibres of 
$\Aut(\IM)$ (each of which is isomorphic to $G$),

$\bullet$\ Similarly it is a {\em Levi reduction} if the fibres of $\cP$ 
embed as Levi subgroups of parabolic subgroups 
of the fibres of $\Aut(\IM)$,

$\bullet$\ The reduction is {\em proper} if the fibres
of $\cP$  embed as proper subgroups.
\end{defn}

\begin{lem}\label{lem: reducible and ppr}
$\IL$ is reducible if and only if it has a proper parabolic reduction of structure group.
\end{lem}
\pf
Given $\cP\subset \Aut(\IL)$ with $\IT\into \cP\bigl\vert_\IH$,
then taking $\IP=\cP$ and $\IM=\IL$ gives the desired reduction.
Conversely given $\IM,\cP,\IP$ then  $\Aut(\IP)$
gives the desired parabolic
sublocal system in $\Aut(\IL)$. 
\epf

Note that part 1) of Thm. \ref{thm: red-irred-SLS}
now follows immediately from Cor. \ref{cor: irredGalrho}, noting that 
$\Gal(\IL)=\Gal(\rho)$.

This irreducibility condition can also be spelt out in terms of Stokes representations and compatible systems of parabolics 
(as in \cite{gbs} \S9). %

If $\cG$ is constant then we can use the usual (simpler) 
notion of reduction of structure group 
(then $\Gal\subset G$ 
and we don't need twisted reductions, 
i.e. we can take $\IM$ to be the trivial $G$-torsor).

\subsection{Reductive/semisimple graded local systems}\label{ss: red-ss loc sys}

If $\IL\to S$ is a graded local system and $\cL\to S$ 
is a Levi reduction of $\IL$ 
(in the sense of Defn. \ref{defn: twisted reductions of structure}) 
then $\cL$ is itself a graded local system, 
and so we can ask if $\cL$ is reducible or not.
Define $\IL$ to be {\em reductive} (or ``semisimple'') 
if it has an irreducible Levi reduction.
Similarly to Lemma 
\ref{lem: reducible and ppr}
one can rephrase this in terms of $\Aut(\IL)$:

\begin{lem}\label{lem: reductive and plr}
$\IL$ is reductive 
if and only if there is a sublocal system $\cE\subset \Aut(\IL)$ 
such that 1) each fibre $\cE_p$ is a Levi subgroup of a parabolic of $\Aut(\IL)_p$, 2) $\IT\subset \cE$, and 3) $\cE$ is irreducible in the sense that it has no proper sublocal systems of parabolic subgroups, containing $\IT$.
\end{lem}

Recall $\bar\Gal(\IL)\subset \Ga\subset \Aut(G)$ is the image of 
$\Gal(\IL)$ in $\Ga$.

\begin{prop}\label{lem: 2red-sls-defns}
$\IL$ is reductive if and only if $\bar\Gal(\IL)$ is a linearly reductive group.
\end{prop}
\pf
By Prop. \ref{prop: Lred-levi} $\bar\Gal(\IL)$ is a linearly reductive group if and only if there is a subgroup $L\subset G$ such that 1) $L$ is a Levi subgroup of a parabolic subgroup of $G$, 2) $L$ is preserved by the action of $\Gal(\IL)$, and 3) $L$ has no proper parabolic subgroups that are $\Gal(\IL)$ invariant. 
(This uses the fact that the centralisers of tori in $G$ are exactly the Levi subgroups of parabolics.)
As in Lem. \ref{lem: 2irred-sls-defns}  the existence of such $L$ is the same as $\IL$ having an irreducible Levi reduction.
\epf

Part 2) of Thm. \ref{thm: red-irred-SLS}
now follows immediately from Cor. \ref{cor: main-red-annonce}.

Note that if $\cG$ is constant (or has finite monodromy) 
then this is the same as $\Gal(\IL)$ being linearly reductive.

\begin{rmk}
Note that if $\IL$ is a Stokes $\cG$-local system then it makes no difference 
if we insist on only looking at reductions that are Stokes local systems: 
i.e. $\IL$ is semisimple if it has a Levi reduction to an irreducible Stokes 
$\cL$-local system, and  $\IL$ is reducible  
if it has a parabolic reduction to a Stokes $\cP$-local system.
To see this we need to check 1) that the Stokes conditions 
on the monodromies (of the reductions) 
around the tangential punctures are automatic, and 
2) that there is no loss of generality in assuming that the local systems of 
parabolic/Levi subgroups are untwisted around the tangential punctures.
This is now an easy exercise: 
1) follows since the Stokes groups are  controlled by $\IT$, and 2) follows by considering the proof of Lem. \ref{lem: reducible and ppr}, and the analogous proof of Lem. \ref{lem: reductive and plr}. 
\end{rmk}

\begin{rmk}
In the case where $\cG$ is a constant general linear group 
with fibre $G=\GL_n(\IC)$ for some $n$,
then a Stokes $\cG$-local system $\IL$ is equivalent to a Stokes local system $\IV$ of rank $n$ vector spaces, as 
in \cite{tops}.
Then $\IL$ is irreducible if and only if 
$\IV$ has no nontrival proper Stokes sublocal systems, and it is reductive if and only if $\IV=\bigoplus \IV_i$ is the direct sum of some irreducible Stokes local systems $\IV_i$.
\end{rmk}

Note that there are thus many simple criteria 
to ensure points of $\cR$ 
are stable, 
and they will be studied systematically elsewhere. 
For example in the constant $\GL_n(\IC)$ setting,  
all the Stokes local systems on $\bSi$ 
are irreducible if at one of the punctures
the irregular class  
just has one Stokes circle $I\subset \cI$ 
with $\Ram(I)=n$ 
(e.g. if $\slope(I)=k/n$ with $(k,n)=1$ this is 
Katz's irreducibility criterion \cite{katz-dgps} 
(2.2.8)).

\begin{eg}
As a very simple  explicit example consider the wild Riemann surface 
$\bSi=(\IP^1,\infty,\<2x^{3/2}\>)$ underlying the 
Airy equation.
Looking at the corresponding Stokes diagram 
(\cite{stokes1857} p.116,  
or e.g. \cite{twcv} p.1 or the app \cite{ssd}), 
this has $3$ singular directions, so
$\wt \Si\simeq \IP^1\setminus\{4 \text{ points}\}$
and so  
the wild surface group $\Pi=\Pi_1(\bSi)$ 
is the free nonabelian group with three generators 
(using just one tangential basepoint at $x=\infty$), and
the wild monodromy relation can be written in the form:
\beq\label{eq: airy relation}
S_3S_2S_1 = h
\eeq
in the group $G=\SL_2(\IC)$, 
where the formal monodromy $h$ and the 
Stokes matrices $S_i$ have the form
$$h=\bmx 0& a \\ -a^{-1} & 0\emx,\qquad
S_1=\bmx  1 & s_1 \\ 0 & 1 \emx,\ 
S_2=\bmx  1 & 0 \\ s_2 & 1 \emx,\ 
S_3=\bmx  1 & s_3 \\ 0 & 1 \emx$$
with $a\in \IC^*, s_i\in \IC$.
The  relation \eqref{eq: airy relation}
then says $s_1=s_3=a, s_2=-a^{-1}$,
so that the space of solutions, 
the wild representation variety 
$\cR=\Hom_\IS(\Pi,G)$,
is a copy of $\IC^*$.
Note that 
this relation \eqref{eq: airy relation}
is one of the central equations in Lie theory\footnote{ 
defining Tits's famous trijection 
\cite{titsstrconstants} p.25, 
and about which 
Faltings wrote 
``... everything comes down to this famous all important formula''
\cite{faltingsloop} p.53 
(although usually it is not recognised as a 
wild monodromy relation).}.
The fact that the diagonal torus acts transitively (by diagonal conjugation) on $\cR$
is the rigidity of the Airy equation 
(i.e. the wild character variety is a point).
For any $a$, 
the group generated by $S_1,S_2,S_3,h$ is dense in
$\SL_2(\IC)$ so Ramis tells us this is the  
Galois group, so the Airy equation is irreducible.
In general we would also need to know the Ramis (exponential) torus as well, so lets compute that in this example
(the notation we use is in \cite{twmcg} for example):
If $b$ is the tangential basepoint at $\infty$ then the 
fibre over $b$ of the Stokes circle 
$\<2x^{3/2}\>$ has two points, given by the 
two functions $\pm 2x^{3/2}$ and they generate 
a rank one sub-lattice (free $\IZ$ module)
$$X^*(\IT_b) \subset \cI_b$$
in the complex vector space $\cI_b$ of all the 
Fabry functions in the direction $b$.
The Ramis torus $\IT_b\cong \IC^*$ is the torus with 
this character lattice $X^*(\IT_b)\cong \IZ$.
The irregular class further specifies an embedding 
$\IT_b\into \SL_2(\IC)$ (from the grading of the fibre of the graded local system) so we can add $\IT_b$ to the generators of the Galois group.

\end{eg}

\renewcommand{\baselinestretch}{1}              %
\normalsize
\bibliographystyle{amsplain}    \label{biby}
\bibliography{../thesis/syr}

\noindent
Universit\'e Paris Cit\'e and Sorbonne Universit\'e, CNRS,
IMJ-PRG, 75013 Paris, France.\hfill 
boalch@imj-prg.fr\hfill \!
\noindent
Department of Mathematics,
Faculty of Science Division I,
Tokyo University of Science, 
1-3 Kagurazaka, Shinjuku-ku,
Tokyo 162-8601, Japan.\ 
yamakawa@rs.tus.ac.jp

\end{document}